\documentclass[a4paper]{article}
\usepackage[T1]{fontenc}
\usepackage[utf8]{inputenc}	
\usepackage[tbtags]{amsmath}
\usepackage{amsfonts,amssymb,mathrsfs,amscd,comment}

\usepackage{xcolor}
\usepackage{hyperref}
\usepackage[capitalise]{cleveref}
\hypersetup{hypertexnames = false,
	bookmarksdepth = 2, bookmarksopen = true, bookmarksnumbered, colorlinks, 
	linkcolor={red!50!black},
	citecolor={blue!50!black},
	urlcolor={blue!80!black},
	pdfstartview={XYZ null null 1}}

\usepackage{umnbib-e} 

\usepackage{geometry}
\def \A {\mathbb{A}}
\def \L {\mathbb{L}}
\def \P {\mathbb{P}}
\def \Q {\mathbb{Q}}
\def \Z {\mathbb{Z}}
\DeclareMathOperator{\Gal}{Gal}
\DeclareMathOperator{\Pic}{Pic}
\DeclareMathOperator{\Prj}{Proj}
\DeclareMathOperator{\Ker}{Ker}

\begin{document}

\title{\bf\large\MakeUppercase{%
On pairs, triples and quadruples of points on a cubic surface
}}

\author{Sergey Galkin, Pavel Popov}
\date{}

\maketitle

\makeatletter
\makeatother
\footnotetext{The authors are partially supported by \href{https://ms.hse.ru/en}{Laboratory of Mirror Symmetry NRU HSE}, RF Government grant, ag. N 14.641.31.0001}

This note considers the problem of stable birational classification
for products of symmetric powers of a fixed algebraic variety $Y$
irreducible over a field $k$.
Varieties $Y, Y'$ are said to be stably birationally equivalent (over $k$)
if there are numbers $n, n'$ and an isomorphism of function fields $k(Y)(z_1,\ldots,z_n) \simeq k(Y')(z'_1,\ldots,z'_{n'})$.
Let $\Sigma_n$ be a permutation group of an~$n$-element set,
$Y^{(n)} := Y^n/\Sigma_n$ --- $n$-th symmetric power of $Y$,
$Y^{[n]}$ --- Hilbert scheme of~$n$~points on $Y$
(if $Y$ is a smooth surface, then $Y^{[n]}\to Y^{(n)}$ is a resolution of singularities by Fogarty's theorem),
$X\subset\P^3_k$ --- geometrically integral cubic surface over $k$.

{\bf Theorem.} The varieties $X\times X^{(3)}$ and $X\times X^{(4)}$ are stably birationally equivalent over~$k$.

{\bf Lemma 1.} The varieties $X^{(2)}\times X^{(4)}$ and $X^{(2)}\times X^{(3)}\times\P^2$ are birationally equivalent~over~$k$.

{\bf Lemma 2.} The varieties $X^{(2)}$ and $X\times\P^2$ are birationally equivalent over $k$.

Note that there are $k$ and $X/k$, such that $X^{(4)}$ are not stably birationally equivalent~to~$X^{(3)}$.

{\bf Example 1.} The surface $F = \Prj k[X_0,X_1,X_2,X_3]/(X_0^3 + 7 X_1^3 + 7^2 X_2^3 - 2 X_3^3)$.
According to \cite{CM04}, $F$ is smooth over $\Q$ and $F(K)=\emptyset$ for field extensions $K/\Q_7$ of degree 1, 2, 4.
So $F^{(4)}(\Q_7)=\emptyset$. But $F^{(3)}$ has a point over any field.
By theorem of Nishimura~\cite{Nis55} the varieties $F^{(4)}$ and $F^{(3)}$ are not stably birationally equivalent over $\Q_7$.

{\bf Example 2.} Brauer--Severi surfaces ($Y/k$ such that $Y_{\bar{k}} \simeq \P^2_{\bar{k}}$).
According to \cite{Kol16} stable birational classification is described by a boolean function: 
``Is $Y(k)$ an empty set?''.
Brauer--Severi surfaces are birationally equivalent to cubic surfaces with $\Gal_k$-invariant \emph{sixes} ($6$ skew lines).

{\bf Example 3.} Cubic curves.
Let $C\subset\P^2_k$ be a smooth cubic curve with the Jacobian $E:=\Pic^0 C$.
Using the Abel map we see that $C^{(n)}$ is stably birationally equivalent to~$\Pic^n C$,
besides $\Pic^{3n\pm1} C \simeq C$ and $\Pic^{3n} C \simeq E$.
If $C(k)=\emptyset$, then $C^{(3)}$ and $C^{(4)}$ are not stably birationally equivalent.
But there is an isomorphism from $C\times C$ to  $C\times E$ : $(P,Q)\mapsto (P,\mathcal{O}(P-Q))$.
Hence $C\times C^{(4)}$ is stably birationally equivalent to $C\times C^{(3)}$.
Since a cone is stably birationally equivalent to its base,
all these results about cubic curves apply to conical cubic surfaces. 

{\bf Twisted cubic curves.}
We say that six points are \emph{in general position} if no four of them lie in the same plane.
Through six points in general position one can pass a unique twisted cubic curve, 
rational normal curve of degree $3$ (see for example Theorem 1.18 in \cite{Ha92ru}).
Note that any six points on a twisted cubic curve are in general position.
The set $M_3$ of all twisted cubic curves is a $12$-dimensional connected smooth variety, which is isomorphic to~$PGL(4)/PGL(2)$.
We say that a curve is \emph{transversal} to a surface, if their scheme-theoretic intersection is a reduced Artin scheme. 
For any reduced cubic surface~$X\subset\P^3$ 
the subset~$U\subset~M_3$ of twisted cubic curves transversal to $X$ is non-empty and open in Zariski topology.
Let $\tilde{U} \to U$ be a normal \'etale covering with deck transformation group $\Sigma_9$:
a point in $\tilde{U}$ is a twisted cubic curve $T$ transversal to $X$ plus a total order on the set $T\cap X$.
Consider the standard embedding of $\Sigma_2\times\Sigma_3\times\Sigma_4$ in $\Sigma_9$
and denote by~$U_{2,3,4}$ the intermediate quotient $\tilde{U}/(\Sigma_2\times\Sigma_3\times\Sigma_4)$.
A point $(T;B,C,D)\in U_{2,3,4}$ is a twisted cubic~$T$ transversal to $X$
plus a partition of the set $T\cap X$ into three subsets $B=\{B_1,B_2\}$, $C=\{C_1,C_2,C_3\}$, $D=\{D_1,D_2,D_3,D_4\}$.
Consider a map $\pi_{2,4}\colon U_{2,3,4}\to X^{(2)}\times X^{(4)}$ such that 
$\pi_{2,4}(T;B,C,D) = (B,D)$.
If a point $(B,D)$ lies in the image of $\pi_{2,4}$ then the set of six points $B\cup D$ lies on some twisted cubic curve $T$,
and so these points are in general position and such $T$ is unique.
Thus $C = (T\cap X)\backslash(B\cup D)$ and $\pi_{2,4}$ is an open embedding.
Fix a plane $\Pi\subset\P^3$ and two points $O,O'\in\P^3(k)\backslash\Pi$.
Define the subset $V\subset U_{2,3,4}$ by the condition that points $O,O',B_1,B_2$ are in general position.
Consider $W = \{(T;A,B,C,D) | (T;B,C,D)\in V, A\in(T\cap\overline{OB_1B_2})\backslash(B\cup C)\}$.
Define a map $\pi_{2,3,1'}\colon W\to X^{(2)}\times X^{(3)}\times\Pi$ as $\pi_{2,3,1'}(T;A,B,C,D) = (B,C,A')$, where $A' := \overline{O'A}\cap\Pi$.
If $(B,C,A')$ lies in the image of $\pi_{2,3,1'}$, then the set of six points $B\cup C\cup A$ (where $A=\overline{O'A'}\cap\overline{OB_1B_2}$)
lies on a twisted cubic curve $T$, so these points are in general position and such $T$ is unique. 
This implies that $D = (T\cap X)\backslash(B\cup C)$ and $\pi_{2,3,1'}$ is an open embedding.
We see that both $X^{(2)}\times X^{(4)}$ and $X^{(2)}\times X^{(3)}\times\Pi$ have Zariski-open subsets isomorphic to $W$.
Lemma 1 is proved.
The proof of Lemma 2 is similar, using lines instead of twisted cubic curves (see also \cite{Popov18}).
The Theorem is an immediate consequence of the Lemmata.

{\bf Corollary.}
Denote by $\Z[SB_k]$ the free abelian group generated by the classes $\{Y\}$ of stable birational equivalence
of varieties $Y$ irreducible over $k$. It can be equipped with operations of multiplication and symmetric powers.
The Theorem implies that for all $X$ we have an identity
$\{X\}\cdot\{X^{(4)}\}=\{X\}\cdot\{X^{(3)}\}\in\Z[SB_k]$, but the Examples show that $\{X^{(4)}\}\neq\{X^{(3)}\}\in\Z[SB_k]$ for some $X$.
These cubic surfaces provide examples of zero divisors in the ring $\Z[SB_k]$.
The Grothendieck ring of $k$-varieties $K(Vars/k)$ is generated by the classes $[Y]$ of schemes $Y$ of finite type over $k$ and relations 
$[Y] = [W] + [Y\backslash W]$ for closed subschemes $W\subset Y$. 
Denote by $\L := [\A^1]$ the class of an affine line. 
Define the quotient-ring $LL(k) := K(Vars/k)/(\L)$.
According to \cite{LL03} in characteristic $0$ there is an isomorphism 
$\mu:LL(k)\to\Z[SB_k]$, $\mu([Y])=\{Y\}$ for all smooth connected projective schemes $Y/k$.
Surfaces from Example 1 and 2 are also zero-divisors in the ring $LL(k)$.
In \cite{Popov18} it is proven that there is a unique relation of degree $4$ between $\Gal_k$-representations $H(X^{[n]})$ and $H(M_3(X))$,
where $M_3(X)$ is a variety of generalized twisted cubic curves on a smooth $X$. 
This relations is not identically true in $K(Vars/k)$, since modulo $(\L)$ it becomes 
$[X^{(4)}] = [X^{(3)}] \in LL(k)$, and this contradicts to Examples 1 and 2.
We expect that a modification of the relation from \cite{Popov18}
by multiplication on $X^{[2]}$ might be a true identity in $K(Vars/k)$.
Theorem implies that in characteristic $0$ for smooth $X$ there is a pair of smooth $10$-dimensional manifolds $N_{\pm}$, such that 
$[X^{[2]}]\cdot([X^{[4]}] - [X^{[3]}]\cdot[\P^2]) = \L\cdot([N_+]-[N_-]) \in K(Vars/k)$.
To find $N_\pm$ one need to decompose the birational morphism $\pi_{2,4}\circ{(\pi_{2,3,1'})}^{-1}$
into a sequence of blow-ups and contractions with smooth centers and describe these centers.

We thank K. Shramov for the reference \cite{Nis55},
and C. Voisin, C. Peskine, F. Han for their interest in our work.

\providecommand{\arxiv}[1]{\href{http://arxiv.org/abs/#1}{\tt arXiv:#1}}

\bigskip

{\bf \href{https://mccme.ru/~galkin}{Sergey Galkin}}

\href{https://hse.ru/en/staff/sergey}{HSE University}, Russian Federation

\href{https://ium.mccme.ru/english}{Independent University of Moscow}

{\bf \href{https://www.hse.ru/en/org/persons/14288573}{Pavel Popov}}

\href{https://math.hse.ru/en}{HSE University}, Russian Federation

\href{https://ium.mccme.ru/english}{Independent University of Moscow}

{\it E-mail}: lightmor@gmail.com


\def \R {\mathbb{R}}
\def \C {\mathbb{C}}
\def \CC {\mathcal{C}}
\def \II {\mathcal{I}}

\bigskip

{\bf Afterword} for \arxiv{1810.07001} version of the article.
The volume of the original note is bounded by the journal guidelines,
so we picked but one most surprising fact.
The reader might find useful the commentary below.

We thank J.-L.\,Colliot-Th\'el\`ene for pointing out our attention to the subtlety of defining
$\Z[SB_k]$: if $k$ is algebraically closed one can introduce operations
of multiplication and symmetric powers directly on the set $SB_k$ of stably birational equivalence classes,
however over a non-algebraically closed field $k$ it happens that (symmetric) products of (non-geometrically)
irreducible varieties may fail to be irreducible (think of $\C\otimes_\R\C$),
so product operations applied to basic elements might produce non-trivial linear combinations in $\Z[SB_k]$.
By abuse of notation $\Z[SB_k]$ is often called a ``semi-group ring'' despite $SB_k$ not being an honest semi-group.
Equivalently one can define $\Z[SB_k]$ as a quotient of a free abelian group generated by classes $\{Y\}$
of all (non necessarily irreducible over $k$) varieties $Y/k$ by relations
$\{Y\coprod Y'\}=\{Y\}+\{Y'\}$ for all $Y,Y'/k$ and $\{Z\}=\{Z'\}$ for all $Z,Z'/k$ that are irreducible over $k$
and stably birationally equivalent over $k$; one then can conveniently use classes $\{Y\}$
of reducible varieties.

The note is constrained with geometrically irreducible cubic surfaces for the similar reason:
traditionally one speaks about (stable) birational equivalence of irreducible varieties,
however symmetric cube of a planes triple is reducible; the construction in the proof still
gives (stable) birational equivalences between some irreducible components.
Two distinct alternatives to the notion of stable birational equivalence are the equalities of the classes in 
either $\Z[SB_k]$ or $LL(k)$.
To deal with non-reduced cubic surfaces one has to use nested Hilbert schemes instead
of the products (for smooth surfaces these varieties are birationally equivalent).

Consider a polynomial ring $R := \Z[h_1,h_2,\dots]$.
Every variety $Y$ over a field $k$ gives rise to two ring homomorphisms
$\phi_{Y,k} : R\to\Z[SB_k], \psi_{Y,k} : R \to LL(k)$ defined by
$\phi_Y(h_i) = \{Y^{(i)}\} \in \Z[SB_k], \psi_Y(h_i) = [Y^{(i)}] \in LL(k)$.
The two homomorphisms agree in characteristic $0$ if $Y$ is smooth proper,
however they disagree in Example 3.
For a class $\CC$ of fields $k$ and varieties $Y/k$ one may consider the homomorphic images
of $\phi_\CC := \prod_{(Y/k)\in\CC} \phi_{Y,k}$ and $\psi_\CC := \prod_{(Y/k)\in\CC} \psi_{Y,k}$.
The search for ``beautiful formulae'' in particular contributes to the study of these images
and of the respective kernels 
$\II_\phi(\CC) := \bigcap_{(Y/k)\in\CC} \Ker \phi_{Y,k}$,
$\II_\psi(\CC) := \bigcap_{(Y/k)\in\CC} \Ker \psi_{Y,k}$.
Main results immediately imply that $h_1(h_4-h_3)$ belongs to
$\II_\phi(GICS)\cap\II_\psi(SCS)$ where $SCS$ (resp. $GICS$) is the class of all
smooth (resp. geometrically integral) cubic surfaces,
however extra analysis is needed to see whether it also belongs to $\II_\psi(GICS)$.
One can introduce equivalence relation on natural numbers $n$ and partitions $\lambda=(n_1,\dots,n_k)$
by saying that $\lambda \sim \lambda'$ iff $\phi_\CC(s_\lambda - s_{\lambda'})=0$,
here $s_\lambda := \prod h_i^{n_i}$.\footnote{We have to note that all 
pre-$\lambda$ rings $K_0(Vars/k), LL(k), \Z[SB_k]$ fail to be special $\lambda$-rings,
so in the stable birational setting some terminology of $\lambda$-rings could be a false friend.}
Main theorem is the equivalence $(4,1)\sim(3,1)$,
Lemma 2 implies $2\sim 1$,
a variation of the construction with twisted cubics also shows $4\sim 5$ and $3\sim 6$.
Similar considerations of space curves of genus one and degree four (base loci for pencils of quadrics)
show that $4\sim 8$ and $5\sim 7$.
So $1\sim 2$, $4\sim 5\sim 7\sim 8$ and $3\sim 6$, but $3k$ is never equivalent to $3l\pm1$.
J.\,Koll\'ar kindly pointed to us that $3\sim 6$, and asked whether $3$ is equivalent to $0$ or not.
Example $3$ shows that $3$ is not equivalent to $0$ in the class of geometrically integral cubic surfaces,
so if the claim is true for the class of smooth cubics, then the proof shall at some point explicitly
discriminate the cones.

The proofs of Lemmata in fact imply that $X^{(4)}\times X$ is birationally equivalent to
$X^{(3)}\times X\times\P^2$,
because birational isomorphism from the proof of Lemma 1 commutes with the projection to $X^{(2)}$,
so we can restrict it to a general rational section $X \dashrightarrow X^{(2)}$
of the birational isomorphism $X^{(2)} \dashrightarrow X\times\P^2$ from the proof of Lemma 2.

Similar constructions with curves of degree $3$ and $4$ can be used to equip
$X^{(3)}$ and $X^{(4)}$ with structures of ``birational symmetric quasi-groups''
similarly to Manin's observation that lines equip $X$ itself with a structure of ``birational Moufang loop''.
This produces the families of birational involutions on $Z\in\{X,X^{(3)},X^{(4)}\}$
birationally parameterised by $Z$ itself.

\bigskip

These matters will be addressed in the forthcoming PhD thesis of the second author.


\end{document}


\title{\bf\large\MakeUppercase{%
О парах, тройках и четвёрках точек на кубической поверхности
}}

\author{С.\,С.~Галкин, П.\,П.~Попов}
\date{}

\maketitle

\makeatletter
\renewcommand{\@makefnmark}{}
\makeatother
\footnotetext{Работа выполнена при поддержке Лаборатории зеркальной симметрии НИУ ВШЭ, грант Правительства РФ Договор № 14.641.31.0001.}

Заметка посвящена вопросу стабильной бирациональной классификации
произведений симметрических степеней фиксированного алгебраического многообразия $Y,$
неприводимого над полем $k$.
Говорят, что $Y,Y'$ стабильно бирациональны (над $k$),
если существуют числа $n,n'$ и изоморфизм полей функций $k(Y)(z_1,...,z_n) \simeq k(Y')(z'_1,...,z'_{n'})$.
Пусть $\Sigma_n$ это группа перестановок $n$-элементного множества,
$Y^{(n)} := Y^n/\Sigma_n$ --- $n$-ая симметрическая степень $Y$,
$Y^{[n]}$ --- схема Гильберта $n$ точек на $Y$
(если $Y$ это гладкая поверхность, то $Y^{[n]}\to Y^{(n)}$ это разрешение особенностей по теореме Фогарти),
$X\subset\P^3_k$ --- геометрически целая кубическая поверхность над $k$.

{\bf Теорема.} Многообразия $X\times X^{(3)}$ и $X\times X^{(4)}$ стабильно бирациональны над $k$.

{\bf Лемма 1.} Многообразия $X^{(2)}\times X^{(4)}$ и $X^{(2)}\times X^{(3)}\times\P^2$ бирациональны над $k$.

{\bf Лемма 2.} Многообразия $X^{(2)}$ и $X\times\P^2$ бирациональны над $k$.

Заметим, что существуют $k$ и $X/k$, такие что $X^{(4)}$ стабильно не бирационально $X^{(3)}$.

{\bf Пример 1.} Поверхность $F = \Prj k[X_0,X_1,X_2,X_3]/(X_0^3 + 7 X_1^3 + 7^2 X_2^3 - 2 X_3^3)$.
Следуя \cite{CM04}, $F$ гладкая над $\Q$ и $F(K)=\emptyset$ для всех расширений $K/\Q_7$ степеней 1,2,4.
Поэтому $F^{(4)}(\Q_7)=\emptyset$. Но $F^{(3)}$ имеет точки над любым полем.
По теореме Нишимуры \cite{Nis55} многообразия $F^{(4)}$ и $F^{(3)}$ стабильно не бирациональны над $\Q_7$.

{\bf Пример 2.} Поверхности Севери--Брауэра ($Y/k$ такие что $Y_{\bar{k}} \simeq \P^2_{\bar{k}}$).
Согласно \cite{Kol16} стабильная бирациональная классификация описывается булевой функцией
<<пусто ли множество $Y(k)$>>.
Поверхности Севери--Брауэра бирациональны кубическим поверхностям, имеющим $\Gal_k$-инвариантную шестёрку попарно дизьюнктных прямых.

{\bf Пример 3.} Кубические кривые.
Пусть $C\subset\P^2_k$ --- гладкая кубическая кривая с якобианом $E := \Pic^0 C$.
Отображение Абеля показывает, что $C^{(n)}$ стабильно бирационально $\Pic^n C$,
причём $\Pic^{3n\pm1} C \simeq C$ и $\Pic^{3n} C \simeq E$.
Если $C(k)=\emptyset$, то $C^{(3)}$ и $C^{(4)}$ стабильно не бирациональны.
Но $C\times C$ изоморфно отображается в $C\times E$: $(P,Q)\to (P,\mathcal{O}(P-Q))$.
Поэтому $C\times C^{(4)}$ стабильно бирационально $C\times C^{(3)}$.
Поскольку конус стабильно бирационален своему основанию,
все результаты о кубических кривых эквивалентны аналогичным результатам о кубических поверхностях-конусах.

{\bf Скрученные кубики.}
Говорят, что шестёрка точек в пространстве находится в общем положении, если никакие четыре из них не лежат в одной плоскости.
Через такую шестёрку можно провести единственную скрученную кубику (см. например Теорему 1.18 в \cite{Ha92ru}),
рациональную нормальную гладкую кривую степени $3$.
Также любая шестёрка точек на любой скрученной кубике является шестёркой точек в общем положении.
Множество $M_3$ всех скрученных кубик это $12$-мерное связное гладкое многообразие, изоморфное $PGL(4)/PGL(2)$.
Будем говорить, что кривая трансверсальна к поверхности, если их схемное пересечение это приведённая артинова схема.
Для каждой приведённой кубической поверхности $X\subset\P^3$ 
подмножество $U\subset M_3$ трансверсальных к $X$ скрученных кубик непусто и открыто по Зарисскому.
Пусть $\tilde{U} \to U$ это нормальное этальное накрытие с накрывающей группой $\Sigma_9$:
точка $\tilde{U}$ это трансверсальная к $X$ скрученная кубика $T$ плюс полный порядок на множестве $T\cap X$.
Стандартно вложим $\Sigma_2\times\Sigma_3\times\Sigma_4$ в $\Sigma_9$,
и обозначим через $U_{2,3,4} = \tilde{U}/(\Sigma_2\times\Sigma_3\times\Sigma_4)$ промежуточный фактор.
Точка $(T;B,C,D)\in U_{2,3,4}$ это трансверсальная к $X$ скрученная кубика $T$
плюс разбиение множества $T\cap X$ на три группы $B=\{B_1,B_2\}, C=\{C_1,C_2,C_3\}, D=\{D_1,D_2,D_3,D_4\}$.
Определим отображение $\pi_{2,4}:U_{2,3,4}\to X^{(2)}\times X^{(4)}$ так:
$\pi_{2,4}(T;B,C,D) = (B,D)$.
Если $(B,D)$ лежит в образе $\pi_{2,4}$, то шестёрка точек $B\cup D$ лежит хотя бы на одной скрученной кубике $T$,
а следовательно находится в общем положении и лежит на единственной скрученной кубике;
поэтому $C = (T\cap X)\backslash(B\cup D)$ и $\pi_{2,4}$ это открытое вложение.
Зафиксируем плоскость $\Pi\subset\P^3$ и пару точек $O,O'\in\P^3(k)\backslash\Pi$.
Определим подмножество $V\subset U_{2,3,4}$ условием, что точки $O,O',B_1,B_2$ находятся в общем положении.
Пусть $W = \{(T;A,B,C,D) | (T;B,C,D)\in V, A\in(T\cap\overline{OB_1B_2})\backslash(B\cup C)\}$.
Морфизм $\pi_{2,3,1'}: W\to X^{(2)}\times X^{(3)}\times\Pi$ задан так: $\pi_{2,3,1'}(T;A,B,C,D) = (B,C,A')$, где $A' := \overline{O'A}\cap\Pi$.
Если $(B,C,A')$ лежит в образе $\pi_{2,3,1'}$, то шестёрка точек $B\cup C\cup A$ (где $A = \overline{O'A'}\cap\overline{OB_1B_2}$)
лежит на скрученной кубике $T$, следовательно находится в общем положении и лежит на единственной скрученной кубике,
поэтому $D = (T\cap X)\backslash(B\cup C)$ и $\pi_{2,3,1'}$ это открытое вложение.
Таким образом, и $X^{(2)}\times X^{(4)}$, и $X^{(2)}\times X^{(3)}\times\Pi$ имеют открытые по Зарисскому подмножества, изоморфные $W$.
Лемма 1 доказана.
Лемма 2 доказывается аналогично, с использованием прямых вместо скрученных кубик (также см. \cite{Popov18}).
Теорема является мгновенным следствием комбинации двух Лемм.

{\bf Следствия.}
Определим $\Z[SB_k]$ как свободную абелеву группу, порождённую классами $\{Y\}$ стабильной бирациональной эквивалентности
неприводимых над $k$ многообразий $Y$. Её можно снабдить операциями умножения и симметрических степеней.
Теорема и примеры влекут, что для всех $X$ выполнено тождество
$\{X\}\cdot\{X^{(4)}\}=\{X\}\cdot\{X^{(3)}\}\in\Z[SB_k]$, но существуют $X$ для которых $\{X^{(4)}\}\neq\{X^{(3)}\}\in\Z[SB_k]$.
Такие кубические поверхности доставляют примеры делителей нуля в кольце $\Z[SB_k]$.
Кольцо Гротендика $k$-многообразий $K(Vars/k)$ порождено классами $[Y]$ схем $Y$ конечного типа над $k$ и соотношениями
$[Y] = [W] + [Y\backslash W]$ для замкнутых подсхем $W\subset Y$. 
Пусть $\L := [\A^1]$ это класс аффинной прямой. 
Определим фактор-кольцо $LL(k) := K(Vars/k)/(\L)$.
Согласно \cite{LL03} в характеристике $0$ существует изоморфизм
$\mu:LL(k)\to\Z[SB_k]$, такой что $\mu([Y])=\{Y\}$ для гладких проективных связных $Y/k$.
Поверхности из примеров 1 и 2 также являются делителями нуля в кольце $LL(k)$.
В \cite{Popov18} доказано, что есть лишь одно соотношение степени $4$ между $\Gal_k$-представлениями $H(X^{[n]})$ и $H(M_3(X))$,
где $M_3(X)$ это многообразие обобщённых скрученных кубик на гладкой $X$.
Это соотношение не выполнено в $K(Vars/k)$, потому что по модулю $(\L)$ оно имеет вид
$[X^{(4)}] = [X^{(3)}] \in LL(k)$, что противоречит примерам 1 и 2. Мы ожидаем, что модификация формулы из \cite{Popov18},
получающаяся домножением на $X^{[2]}$ уже может быть выполнена и в $K(Vars/k)$.
Теорема влечёт, что в характеристике $0$ для гладких $X$ существует пара гладких $10$-мерных многообразий $N_{\pm}$, таких что 
$[X^{[2]}]\cdot([X^{[4]}] - [X^{[3]}\cdot\P^2]) = \L\cdot([N_+]-[N_-]) \in K(Vars/k)$.
Чтобы найти $N_\pm$ нужно разложить бирациональный морфизм $\pi_{2,4}\circ{(\pi_{2,3,1'})}^{-1}$
в последовательность раздутий и стягиваний с гладкими центрами,
и описать эти центры.

Мы благодарим К. Шрамова за ссылку \cite{Nis55},
а также К.Вуазен, К.Пескина, Ф.Хана за проявленный интерес к нашей работе.

\providecommand{\arxiv}[1]{\href{http://arxiv.org/abs/#1}{\tt arXiv:#1}}

\bigskip

{\bf С.\,С.~Галкин (S.\,S.~Galkin)}

Национальный исследовательский университет <<Высшая школа экономики>>

{\bf П.\,П.~Попов (P.\,P.~Popov)}

Национальный исследовательский университет <<Высшая школа экономики>>

{\it E-mail}: lightmor@gmail.com